\documentclass[12pt]{article}
\usepackage[english]{babel}
\usepackage{amsmath,amstext,amsfonts,amssymb,fancyhdr}
\usepackage[T1]{fontenc}

\newcommand{\R}{\mathbb R}

\newcommand\eps{\varepsilon}
\def\A{{\mathcal A}}
\def\B{{\mathcal B}}
\def\F{{\mathcal F}}
\def\HH{{\mathcal H}}
\DeclareMathOperator{\dive}{div}
\DeclareMathOperator{\dist}{dist}

\newtheorem{thm}{Theorem}

\newtheorem{dfn}[thm]{Definition}

\newenvironment{prf}%
{\par\noindent{\it Proof:}\nopagebreak\normalsize}%
{\hfill\linebreak[2]\hspace*{\fill}$\Box$\\[12pt]}

\begin{document}
\title{Stationary Configuration of some Optimal Shaping}
\author{Giuseppe Buttazzo and Al-hassem Nayam}

\date{}
\maketitle

\noindent{\bf Abstract:} We consider the problem of optimal location
of a Dirichlet region in a $d$-dimensional domain $\Omega$ subjected
to a given right-hand side $f$ in order to minimize some given
functional of the configuration. While in the literature the
Dirichlet region is usually taken $d-1$ dimensional, in this shape
optimization problems, we consider two classes of control variables,
namely the class of one dimensional closed connected sets of finite
one dimensional
 Hausdorff measure and the class of sets of points
of finite cardinality, and we give a necessary condition of
optimality.
{\bf Biographical notes}: Giuseppe Buttazzo graduated in Mathematics
at the University of Pisa in 1976, and received in the same year the
Diploma of Scuola Normale Superiore di Pisa, where he won a national
competition as a student. He made the Ph. D. studies from 1976 until
1980 at Scuola Normale di Pisa, where he also got his first
permanent position as a researcher in 1980. He obtained the full
professorship in Mathematical Analysis at University of Ferrara in
1986, and moved in 1990 to University of Pisa, where he presently
works at the Department of Mathematics. His main research interests
include calculus of variations, partial differential equations,
optimization problems, optimal control theory. On these subjects he
wrote more than 150 scientific papers and several books.\\\\
Al-hassem Nayam received his PhD Degree in mathematics at University
of Pisa in 2011. His PhD Thesis mainly deals with the shape
optimization problems of higher codimension. His research interest
lies on asymptotic shapes, necessary condition of optimality and
regularity of optimal shapes. At present he works at International
Centre for Theoretical Physics (ICTP) in Trieste, Italy.

\section{Introduction}\label{intro}

We consider the problem of finding the
optimal location of a Dirichlet
region $\Sigma$ in a $d$-dimensional
domain $\Omega$ associated the $p$-Laplacian
\begin{equation}\label{eq1}
\begin{cases}
-\dive\big(|\nabla u|^{p-2}\nabla u\big)=f&\hbox{in }\Omega\setminus\Sigma\\
u=0&\hbox{on }\partial\Omega\cup\Sigma,
\end{cases}
\end{equation}
where the right hand side $f$ is given as a
nonnegative element of $L^{p'}(\Omega)$, being
$p'$ the conjugate exponent of $p$. The functional
$\F$ we consider as a cost is defined by
$$\F(\Sigma)=\int_\Omega F\big(x,u(x),\nabla u(x)\big)\,dx$$
where $u$ is the unique solution of equation
\eqref{eq1} and $F:\Omega\times\R\times\R^d\to\R$
is a Carath\'eodory function. The shape
optimization problems we consider consists
in the minimization of the
functional $\F$ over two classes of admissible
control variables $\Sigma$. The first class consists
of all closed connected subsets $\Sigma$ of $\Omega$
whose one-dimensional Hausdorff measure $\HH^1(\Sigma)$
(sometime called the total length of $\Sigma$)
is uniformly bounded by some constant
$L$, while the second is the class of discrete
subsets of $\Omega$ whose $0$-dimensional Hausdorff
measure (i.e. their cardinality) does not
exceed a given number $n$.

The existence of an optimal configuration for
 the two optimization problems described
above can be obtained by using a generalization
of the \v{S}ver\'ak's compactness result
(see \cite{Sverak} for $p=2$ and \cite{BucTre}
for a general $p$). In this paper, we consider
a penalized version of these shape optimization
problems \big(see \eqref{eqs1} and \eqref{eqs2}\big)
 by adding to the cost functional $\F$ a Lagrange
multiplier penalization of the form
$\lambda\HH^1(\Sigma)$ and $\lambda\HH^0(\Sigma)$
respectively. This problem has been considered
in \cite{ButMaiSte} in the simplest case when
$p=2$, where the PDE \eqref{eq1} reduces to
the classical Laplace equation. Moreover, it
has been shown (see \cite{busa07}) that, for
the particular case where $F(x,u,\xi)=f(x)u$,
the limit problem, as $p\to+\infty$, is the
minimum problem for the average distance functional
$$\F(\Sigma)=\int_\Omega\dist(x,\Sigma\cup\partial\Omega)f(x)\,dx,$$
where $\dist(x,E)$ denotes the distance
between the point $x$ and the set $E$.

For the well posedness of the minimum problems
\eqref{eqs1} and \eqref{eqs2} the constraint
 $p>d-d_\Sigma$, where $d_\Sigma$ is the dimension
 of $\Sigma$, has to be imposed in order to have
Dirichlet regions of positive $p$-capacity. Let
 us mention that the regularity of the optimal
sets (in the case of problem \eqref{eqs1}) is
still an open problem, even in the simplest two
dimensional setting with $p=2$ and $F=f(x)u$.

\section{Setting of the problem and existence of minimizers}

We consider the following classes of control variables:
$$\A(\Omega):=\big\{\Sigma\subset\Omega\ :\ \Sigma\mbox{ closed connected},\ \HH^1(\Sigma)<+\infty\big\}$$
$$\B(\Omega):=\big\{\Sigma\subset\Omega\ :\ \Sigma\mbox{ discrete},\ \HH^0(\Sigma)=\#(\Sigma)<+\infty\big\}$$
and the shape optimization problems
\begin{equation}\label{eqs1}
\min\big\{\F(\Sigma)+\lambda\HH^1(\Sigma)\ :\ \Sigma\in\A(\Omega)\big\}
\end{equation}
\begin{equation}\label{eqs2}
\min\big\{\F(\Sigma)+\lambda\HH^0(\Sigma)\ :\ \Sigma\in\B(\Omega)\big\}.
\end{equation}
The penalization terms $\lambda\HH^1(\Sigma)$
and $\lambda\HH^0(\Sigma)$ with $\lambda>0$
replace the constraint on $\HH^1(\Sigma)$ and
$\HH^0(\Sigma)$ and prevent the minimizing
sequences to spread over all the domain $\Omega$
and hence getting a trivial solution. The existence
of minimizers in the two shape optimization problems
 $\eqref{eqs1}$ and $\eqref{eqs2}$ is a consequence
of the \v{S}ver\'ak continuity-compactness result
(see \cite{Sverak} for $p=2$ and \cite{BucTre}
for general $p$) and the Blaschke and Go\l ab
theorems. For the convenience of the reader
let us give some details on the existence of an optimal
shape. Let $\{\Sigma_{n}\}_{n}\subset \A(\Omega)$
be a minimizing sequence
in the optimization problem \eqref{eqs1}, then there exists
a constant $C>0$ such that
$\sup_{n}\{\F(\Sigma_{n}) + \HH^{1}(\Sigma_{n})\} \leq C.$
Since $\{\Sigma_{n}\}_{n}$ is a sequence a closed connected
subsets of $\Omega$ such that
$\sup_{n}\HH^{1}(\Sigma_{n}) \leq C$,
by Blaschke theorem (compactness of the sequence
$\{\Sigma_{n}\}_{n}$ in the Hausdorff topology)
and by Go\l ab theorem (lower semicontinuity of the
$\HH^1$ with respect to the Hausdorff topology),
up to extracting a subsequence, $\{\Sigma_n\}_n$
converges in Hausdorff distance to some
$\Sigma\in\A(\Omega)$ and $\HH^1(\Sigma)\le\liminf_{n
\to\infty}\HH^1(\Sigma_n)$. For the
lower semicontinuity of the energy part of the
functional we need the \v{S}ver\'ak
continuity-compactness result which is stated in
this terms: let $\{\Omega_n\}_n$ be a sequence of
open and bounded sets contained in a fix bounded set
$D$. If we assume that the number of the connected
components of $D\setminus\Omega_n$ is
uniformly bounded by some number $k$, then
$\{\Omega_n\}_n$ converges in the Hausdorff
topology to some open and bounded set
$\Omega\subset D$ and the number of the connected
components of $D\setminus\Omega$ is less or equal
to $k$. Moreover, if we denote by $u_{n}
\in W^{1, p}_{0}(\Omega_{n})$ the distributional
solution of the $p$-Laplace equation
$$\begin{cases}
-\Delta_pu_n=f&\hbox{in }\Omega_n\\
u_n\in W^{1,p}_0(\Omega_n)
\end{cases}$$
for some $f$ in $W^{-1, p'}(D)$, then up to subsequence,
$u_{n}$ converges strongly in $W^{1,p}(D)$
($u_{n}$ are extended by zero outside $\Omega_{n}$) to the
function $u$ which is the distributional solution
of the equation
$$\begin{cases}
-\Delta_pu=f&\hbox{in }\Omega\\
u\in W^{1,p}_0(\Omega).
\end{cases}$$
This result is interesting only in the case where $p$
satisfies $d-1 < p \leq d$ because the case where
$p>d$ is trivial due to the fact that functions in
$W^{1, p}(D)$ are continuous and the convergence of
solutions follows easily. To apply this result to
our problem, we choose $\Omega_n=
D\setminus\Sigma_n$ and notice that
$\{\Omega_n\}_n$ converges to $\Omega=D\setminus\Sigma$
in the Hausdorff topology where $\Sigma$ is the limit
of $\Sigma_n$. From the hypothesis
of the function $F$ and the continuity with respect
to the domains variation of solutions, the lower
semicontinuity follows easily
and also the existence of an optimal shape. The
existence of an optimal set in the problem
\eqref{eqs1} is even easier due to the fact that
$p>d$.

Our goal is to derive first order necessary
conditions of optimality, assuming that the
solutions $\Sigma$ of the minimum problems are
regular as necessary. Notice that, since every
set in $\Sigma\in\A(\Omega)$ is countably $\HH^1$
rectifiable, it can be written as an union of
countable $C^1$ curves and a set of $\HH^1$ measure zero.

We distinguish three cases according to
technical computations. We
consider the following classes of control variables
$\Sigma$: the class of closed connected subsets of
$\R^2$, the class of points of $\R^d$ for $d>1$,
and the class of closed connected subsets of
$\R^d$ for $d>2$. In the two last cases, some
extra difficulties occur because of the
co-dimension of $\Sigma$ which is greater
than $1$. For simplicity, we assume that
$\Omega$ has a smooth boundary and $u=0$ on
 $\partial\Omega\cup\Sigma$; we also assume
as much as needed the regularity on the
data. Before looking for the necessary
conditions of optimality, we recall some
definitions and results which will be helpful;
we refer to \cite{BouButFra} for the details.

For a measure $\mu$ we denote for $\mu$ a.e. $x$
by $P_\mu(x,\cdot):\R^d\to\mathrm{Tan}(\mu,x)$
the orthogonal projection of $\R^d$ on $\mathrm{Tan}(\mu,x)$ where
$\mathrm{Tan}(\mu,x)$ stands for the tangent space
of $\mu$ at $x$ that is the set of all tangent
measures to $\mu$ at the point $x$ (see for instance \cite{AmbFusPal}).

\begin{dfn}\label{d1}
The curvature of $\mu$ is defined as
the vector valued distribution
$$H_\mu:=\dive(P_\mu{\mu}).$$
In other words $H_\mu$ is defined by
$$\langle H_\mu,X\rangle=-\int_{\R^d}\dive^\mu X\,d\mu\qquad\forall X\in
C^\infty_c(\R^d),\R^d),$$
where $\dive^\mu X=\sum_{j=1}^d(P_\mu(\nabla X^j))_j$.
\end{dfn}

We denote by ${\cal M}_{BC}$ the set of all
positive and finite Borel regular measures
of $\R^d$ whose curvature is a Borel
regular measure with
finite total mass. Since the curvature
$H_\mu$ of a measure $\mu\in{\cal M}_{BC}$
is not necessary absolutely continuous with
respect to $\mu$, by Radon-Nikodym theorem, we can write
$$H_\mu=h(\mu)\mu+\partial\mu,$$
where $h(\mu)\in L^1_\mu(\R^d,\R^d)$ is the
density of $H_\mu$ with respect to $\mu$
(also called the pointwise curvature) and
$\partial\mu$ is the singular part of
$H_\mu$ with respect to $\mu$ (also called the boundary of $\mu$).

If $\mu=\HH^k\llcorner\Sigma$, with
$\Sigma$ a $C^2$ $k$-manifold with boundary
 in $\R^d$, then by classical divergence
theorem we have
$$H_\mu=\nu\HH^{k-1}\llcorner\partial\Sigma+h\HH^k\llcorner\Sigma,$$
where $h$ stands for the mean curvature
vector of $\Sigma$ and $\nu$ the co-normal
 unit vector of $\partial\Sigma$. When the
 tangent space to $\mu$ is reduced to zero
$\mu$ a.e., $H_\mu$ is zero. This is for
instance the case where $\mu$ is a finite
sum of Dirac masses, or $\mu$ is concentrated
on an $\alpha$-dimensional Cantor subset $C$ of
 $[0,1]$ with $0<\alpha<1$ and $\HH^\alpha(C)\in(0,+\infty)$.

\begin{dfn}\label{d2}
Let $\Sigma$ be a countably $\HH^k$ rectifiable
set and $\mu=\theta\HH^k\llcorner\Sigma$ be the
associated rectifiable measure. A function
$h\in L^1_\mu(\Sigma,\R^d)$ is said to be
the generalized mean curvature of $\Sigma$ if
$$\int_{\R^d}\dive^\Sigma X\,d\mu=-\int_{\R^d}X\cdot h\,d\mu
\qquad\forall X\in C^\infty_c(\R^d,\R^d).$$
\end{dfn}

In this case we denote the generalized
mean curvature of $\Sigma$ by $H_\Sigma$.

\begin{thm}\label{t2}
Let $(\mu_r)_r$ be a bounded sequence in
${\cal M}_{BC}$ weakly converging to some
measure $\mu$ and assume that
$\mathrm{dim Tan}(\mu_r)\mu_r$ weakly
 converges to $g\mu$. Then the condition
\begin{equation}\label{eq2}
\mathrm{dim Tan}(\mu,x)\le g(x)\qquad\mu-a.e.
\end{equation}
is necessary and sufficient to have
\begin{equation}\label{eq3}
P_{\mu_r}\mu_r\rightharpoonup P_\mu\mu.
\end{equation}
In this case we have
\begin{equation}\label{eq4}
H_{\mu_r}\rightharpoonup H_\mu.
\end{equation}
\end{thm}

\begin{prf}
see \cite{BouButFra}
\end{prf}

Let $\Omega$ be an open subset of $\R^d$
and $F:\Omega\times\R\times\R^d\to[0,+\infty]$
be a positive Carath\'eodory function. We
assume $F$ smooth and satisfying the growth condition
$$F(x,u,z)\le a(x)+|u|^p+|z|^p,$$
where $a$ is an $L^1(\Omega)$ function. We
consider the functional
$$\F(\Sigma):=\int_{\Omega\setminus\Sigma}
F\big(x,u(x),\nabla u(x)\big)\,dx+\lambda m(\Sigma)$$
where $m(\Sigma)$ is either $\HH^1(\Sigma)$
if $\Sigma$ is a closed connected one
dimensional set or $\HH^0(\Sigma)$ if $\Sigma$
is a discrete set of points and $u$ the
solution of the equation \eqref{eq1}. We
are interested in the necessary conditions
of optimality satisfied by the minimizers
of $\F$. From now on we assume optimal sets
in the case of closed connected sets to be
of class $C^{1,\alpha}$ for some $0<\alpha\le1$
that is locally graph of $C^{1,\alpha}$ functions.

\section{Case of closed connected subsets in $\R^2$}

Let $u$ be the weak solution of the state equation
\begin{equation}\label{eq5}
\begin{cases}
-\dive\big(|\nabla u|^{p-2}\nabla u\big)=f&\hbox{in }\Omega\setminus\Sigma\\
u=0&\hbox{on }\partial\Omega\cup\Sigma,
\end{cases}
\end{equation}
that means in its weak formulation
$$\int_\Omega|\nabla u|^{p-2}\nabla u\nabla v\,dx
=\int_\Omega fv\,dx\qquad\forall v\in W^{1,p}_0(\Omega\setminus\Sigma).$$
We introduce the family of diffeomorphisms
$\varphi_\eps(x)=x+\eps X(x)$ where $X$ is
a smooth vector field from $\R^d$ to $\R^d$
supported in $\Omega$. For $\eps$ small enough,
 $\varphi_\eps$ maps $\Omega$ into $\Omega$. Set
 $\Sigma_\eps=\varphi_\eps(\Sigma)$ and
consider the new state equation in the
deformed domain
\begin{equation}\label{eq6}
\begin{cases}
-\dive\big(|\nabla u_\eps|^{p-2}\nabla u_\eps\big)=f&\hbox{in }\Omega\setminus\Sigma_\eps\\
u=0&\hbox{on }\partial\Omega\cup\Sigma_\eps.
\end{cases}
\end{equation}
The corresponding functional is
\begin{equation}\label{eq7}
\F(\Sigma_\eps)=\int_{\Omega\setminus\Sigma_\eps}
F\big(x,u_\eps(x),\nabla u_\eps(x)\big)\,dx+\lambda\HH^1(\Sigma_\eps).
\end{equation}
Notice that for $\eps=0$ equation \eqref{eq6}
 reduces to \eqref{eq5} and we denote $u_0$
simply by $u$. From now on $u$ stands for the
unique solution of \eqref{eq5} and $u_\eps$
for \eqref{eq6}. To differentiate \eqref{eq7}
 we need to show the differentiability of the
 function $\eps\mapsto u_\eps$ at zero. If we
assume $f\in L^\infty(\Omega)$ then, by the
regularity theory of elliptic equations, $u$
and $u_\eps$ are in $W^{2,p}_0(\Omega\setminus\Sigma)$
and $W^{2,p}_0(\Omega\setminus\Sigma_\eps)$
respectively. Now let us denote by
 $U_\eps=u_\eps\circ\varphi_\eps$ the
transported solution on the fixed domain
$\Omega\setminus\Sigma$. Since
$U_\eps\in W^{2,p}(\Omega\setminus\Sigma)$
(because $\varphi_\eps$ is smooth and
 $u_\eps\in W^{2,p}(\Omega\setminus\Sigma_\eps)$)
the Lemma 4.2 of \cite{MurSim} gives the
differentiability of $\eps\mapsto U_\eps$
from $[0,\eps_0)$ to $W^{2,p}(\Omega\setminus\Sigma)$
at zero and the differentiability of
 $\eps\mapsto u_\eps$ from $[0,\eps_0)$
to $W^{1,p}_{loc}(\Omega\setminus\Sigma)$
is obtained via Lemma 2.1 of \cite{Simon}. Moreover
$u'=U'-\nabla u\cdot X$ in $\Omega\setminus\Sigma$,
where $u'=\lim_{\eps\to0}(u_\eps-u)/\eps$ in
the distributional sense. For an open set
$D$, taking $v,\varphi\in W^{1,p}(D)$ and
$\psi\in C^\infty_c(D)$, an easy
computation of the limit of
$$\frac{1}{t}\langle\Delta_p(v+t\varphi)-\Delta_p(v),\psi\rangle,$$
as $t\to0$, gives $\langle-\dive G_v(\nabla\varphi),\psi\rangle$ where $\langle\cdot,\cdot\rangle$ is the duality pairing between ${\cal D'}(D)$ and $C^\infty_c(D)$, and $G_v$ is defined by
$$G_v(Z)=|\nabla v|^{p-2}Z+(p-2)|\nabla v|^{p-4}(\nabla v\cdot Z)\nabla v.$$
So the differential of $\Delta_p$ at a point $v\in W^{1,p}(D)$ is given by
$$\frac{\partial\Delta_p(v)}{\partial v}(\varphi)=
-\dive(G_v(\nabla\varphi)),\qquad\forall\varphi\in W^{1,p}(D).$$
At $u$ solution of \eqref{eq5} the differential of $\Delta_p$ is linear and continuous from $W^{1,p}(\Omega\setminus\Sigma)$ into ${\cal D'}(\Omega\setminus\Sigma)$ that is
$$\frac{\partial\Delta_p(u)}{\partial v}\in{\cal L}(W^{1,p}(\Omega\setminus\Sigma);{\cal D'}(\Omega\setminus\Sigma)).$$
The fact that $u_\eps\in W^{2,p}(\Omega\setminus\Sigma_\eps)$ and solves
equation \eqref{eq6}, the differentiability of the $p$-Laplacian operator from $W^{1, p}(\Omega)$ to ${\cal D'}(\Omega)$ and the differentiability
of the transported solution $U_\eps=u_\eps\circ\varphi_\eps$ at $\eps=0$ allow to differentiate the equation \eqref{eq6} in the distributional sense at $\eps=0$ (see \cite{Simon}) and we obtain the equation
$$-\dive\big(G_u(\nabla u')\big)=0\qquad\mbox{in }\Omega\setminus\Sigma$$
in the distributional sense. This gives an equation satisfied by $u'$ inside $\Omega\setminus\Sigma$. It remains to find the boundary condition of $u'$. To this aim, we have to differentiate the boundary condition of
$u_\eps$. Due to the particular setting of our problem (the boundary operator is the identity and $u_\eps=0$ on $\partial\Omega\cup\Sigma_\eps$), the fact that $u\in W^{2,p}(\Omega\setminus\Sigma)$ and the differentiability of the transported solution $U_\eps$ at zero imply the
differentiability of the boundary condition (see \cite{Simon}) and we have
$$u'=0\mbox{ on }\partial\Omega\qquad\mbox{and}\qquad
u'=-\nabla u\cdot X\mbox{ on }\Sigma.$$
The fact that $u'$ vanishes on the boundary of $\Omega$ is due to the compact support of the vector field $X$ in $\Omega$. So the equation satisfied by $u'$ in the distributional sense is
\begin{equation}\label{eq8}
\begin{cases}
-\dive\big(G_u(\nabla u')\big)=0&\mbox{in }\Omega\setminus\Sigma\\
u'=0&\mbox{on }\partial\Omega\\
u'=-\nabla u\cdot X&\mbox{on }\Sigma.
\end{cases}
\end{equation}
To complete this part, let us check the differentiability at $\eps=0$ of the cost function. It is well known (see for example \cite{AmbFusPal})that the length functional $\HH^1(\Sigma_\eps)$ is differentiable at $\eps=0$. The only point to check is concerned with the differentiability of the map $\eps\mapsto\int_\Omega F(x,u_\eps,\nabla u_\eps)\,dx$ at $\eps=0$. The smoothness and the growth condition on $F$ imply that the map $v\mapsto F(\cdot,v,\nabla v)$ maps $W^{2,p}(D)$ into $L^1(D)$ and is differentiable from $W^{1,p}(D)$ into ${\cal D'}(D)$ for any open set $D$. Moreover, if $u$ is the solution of equation \eqref{eq5}, $F(\cdot,u,\nabla u)\in W^{1,1}(\Omega\setminus\Sigma)$ and the map $\eps\mapsto F(\cdot,u_\eps,\nabla u_\eps)\circ\varphi_\eps$ is differentiable thanks to the hypothesis on $u_\eps$ and $U_\eps$. Therefore (see Theorem 3.3 of
\cite{Simon}) we have the differentiability of the cost function at $\eps=0$. Summarizing, by taking the derivative of the functional \eqref{eq6} at $\eps=0$ we get
$$\frac{\partial}{\partial\eps}\big|_{\eps=0}\F(\Sigma_\eps)=
\int_\Omega(F_uu'+F_z\cdot\nabla u'+\dive(FX))\,dx+\lambda\frac{\partial}{\partial\eps}\big|_{\eps=0}\HH^1(\Sigma_\eps),$$
where $u'$ is the solution of the equation \eqref{eq8}. The derivative $\frac{\partial}{\partial\eps}\big|_{\eps=0}\HH^1(\Sigma_\eps)$ that appears in the above variation, according to Theorem 7.31 of \cite{AmbFusPal}, gives
$$\frac{\partial}{\partial\eps}\big|_{\eps=0}\HH^1(\Sigma_\eps)=
\int_\Sigma\dive^\Sigma X\,d\HH^1=-\langle H_\Sigma,X\rangle.$$
As $\Sigma$ is countably $\HH^1$-rectifiable, $\dive^\Sigma$ is the projection of the divergence to the approximate tangent line of $\Sigma$ at $\HH^1$-a.e point of $\Sigma$.

Unfortunately, the quantity $\int_\Omega(F_uu'+F_z\cdot\nabla u')\,dx$ is not easily exploitable. To overcome this problem we introduce the adjoint state equation
\begin{equation}\label{eq9}
\begin{cases}
-\dive\big(G_u(\nabla q)\big)=F_u-\dive(F_z)&\mbox{in }\Omega\setminus\Sigma\\
q=0&\mbox{on }\partial\Omega\cup\Sigma,
\end{cases}
\end{equation}
which has to be understood in the distributional sense
$$\int_\Omega(F_u-\dive(F_z))v\,dx+
\int_\Omega\dive\big(G_u(\nabla q)\big)v\,dx=0\qquad\forall v\in C^\infty_c(\Omega\setminus\Sigma).$$
We are not interested in the regularity of the functions $u$ and $q$ in the whole domain $\Omega$ but only near the optimal set $\Sigma$. The functions $u'$ and $q$ are both in $W^{1,p}(\Omega\setminus\Sigma)$. In the variational formulation of the equation \eqref{eq9} if we take $u'$ as a test function, we have
\begin{equation}\label{eq10}
\int_\Omega(F_u-\dive(F_z))u'\,dx+\int_\Omega\dive(G_u(\nabla q))u'\,dx=0.
\end{equation}

Let $\Omega^+$ and $\Omega^-$ be two sets such that $\Omega=\Omega^+\cup\Omega^-$ and $\Sigma\subset\partial\Omega^+\cap\partial\Omega^-$. The sets $\Omega^+$ and $\Omega^-$ may be obtained by connecting $\Sigma$ to the boundary of $\Omega$ by pieces of smooth curves. The assumption made on $\Sigma$, $\partial\Omega$ and $f$ provide sufficient regularity for $u,u',q$ so that the Green formula can be applied to \eqref{eq10}. For the sequel, we use the following notation: $\nabla u^+$ stands for the trace on $\Sigma$ of $\nabla u$ restricted to $\Omega^+$, $\frac{\partial u^+}{\partial\nu}$ for the trace of the respective normal derivative, $F_z^+\cdot\nu=F_z(x,0,\nabla u^+)\cdot\nu$. Similarly $\nabla u^-$ stands for the trace on $\Sigma$ of $\nabla u$ restricted to $\Omega^-$, $\frac{\partial u^-}{\partial\nu}$ for the trace of the respective normal derivative, $F_z^-\cdot\nu=F_z(x,0,\nabla u^-)\cdot\nu$. Recall also that
$$\nabla u^\pm=\frac{\partial u^\pm}{\partial\nu}\nu$$
because $u^\pm=u=0$ on $\Sigma$ (i.e the tangential derivative of $u^\pm$ over $\Sigma$ vanishes). Let us compute separately the terms of the equation \eqref{eq10}. Starting by the first part, we use the regularity of the boundary of $\Omega^+$ to perform the integration by parts:
$$A^+=\int_{\Omega^+}(F_uu'-\dive(F_z)u')\,dx=
\int_{\Omega^+}(F_uu'+F_z\cdot\nabla u')\,dx-
\int_{\partial\Omega^+}u'F_z\cdot\nu\,d\HH^1,$$
where $\nu$ is the outer normal of $\Omega^+$. It is easy to observe that $u'=0$ on $\partial\Omega\cap\partial\Omega^+$ and $\partial\Omega^+=\Sigma\cup(\partial\Omega^+\setminus(\partial\Omega\cup\Sigma))\cup(\partial\Omega\cap\partial\Omega^+)$ so
$$A^+=\int_{\Omega^+}(F_uu'+F_z\cdot\nabla u')\,dx
+\int_\Sigma\frac{\partial u^+}{\partial\nu}(F_z^+\cdot\nu)X\nu\,d\HH^1
-\int_{\partial\Omega^+\setminus(\partial\Omega\cup\Sigma)}\!\!\!u'F_z\cdot\nu\,d\HH^1.$$
Similarly, taking into account the fact that the outer normal of $\Omega^-$ restricted to $\Omega^+\cap\Omega^-$ is $-\nu$, one gets
$$A^-=\int_{\Omega^-}(F_uu'+F_z\cdot\nabla u')\,dx-\int_\Sigma\frac{\partial u^-}{\partial\nu}(F_z^-\cdot\nu)X\nu\,d\HH^1
+\int_{\partial\Omega^-\setminus(\partial\Omega\cup\Sigma)}\!\!\!u'F_z\cdot\nu\,d\HH^1.$$
Combining the previous relations and using the fact that the two sets $\partial\Omega^+\setminus(\partial\Omega\cup\Sigma)$ and $\partial\Omega^-\setminus(\partial\Omega\cup\Sigma)$ coincide, gives ($A=A^++A^-$)
$$A=\int_\Omega(F_uu'+F_z\cdot\nabla u')\,dx+\int_\Sigma\Big(\frac{\partial u^+}{\partial\nu}F_z^+\cdot\nu-\frac{\partial u^-}{\partial\nu}F_z^-\cdot\nu\Big)X\nu\,d\HH^1.$$
For the second term, the integration by parts leads to
$$\begin{aligned}
B^+&=\int_{\Omega^+}\dive G_u(\nabla q)u'\,dx\\
&=-\int_{\Omega^+}G_u(\nabla q)\cdot\nabla u'\,dx
+\int_{\partial\Omega^+}u'G_u(\nabla q)\cdot\nu\,d\HH^1\\
\end{aligned}$$
where $\nu$ is the outer normal of $\Omega^+$ as in the previous case. It is easily seen that
$$\begin{aligned}
B^+=&-\int_{\Omega^+}G_u(\nabla q)\cdot\nabla u'\,dx
-\int_\Sigma\frac{\partial u^+}{\partial\nu}(G_u(\nabla q)^+\cdot\nu)X\nu\,d\HH^1\\
&+\int_{\partial\Omega^+\setminus(\partial\Omega\cup\Sigma)}
u'G_u(\nabla q)\cdot\nu\,d\HH^1.
\end{aligned}$$
Similarly, we have
$$\begin{aligned}
B^-=&-\int_{\Omega^-}G_u(\nabla q)\cdot\nabla u'\,dx
+\int_\Sigma\frac{\partial u^-}{\partial\nu}(G_u(\nabla q)^-\cdot\nu)X\nu\,d\HH^1\\
&-\int_{\partial\Omega^-\setminus(\partial\Omega\cup\Sigma)}u'G_u(\nabla q)\cdot\nu\,d\HH^1.
\end{aligned}$$
Therefore, summing up one obtains ($B=B^++B^-$)
$$\begin{aligned}
B=&-\int_\Omega G_u(\nabla q)\cdot\nabla u'\,dx\\
&+\int_\Sigma\Big(\frac{\partial u^-}{\partial\nu}G_u(\nabla q)^-\cdot\nu- \frac{\partial u^+}{\partial\nu}{G_u(\nabla q)}^+\cdot\nu\Big)X\nu\,d\HH^1.
\end{aligned}$$
By the linearity of the function $G_u$, we get
$$\int_\Omega G_u(\nabla q)\cdot\nabla u'\,dx
=\int_\Omega G_u(\nabla u')\cdot\nabla q\,dx$$
but, by integration by parts allowed by regularity of $\Omega$ and $\Sigma$, it follows that
$$\int_\Omega G_u(\nabla u')\cdot\nabla q\,dx=
-\int_\Omega\dive(G_u(\nabla u'))q\,dx
+\int_{\partial\Omega\cup\Sigma}qG_u(\nabla u')\cdot\nu\,d\HH^1=0$$
because $u'$ is the weak solution of equation \eqref{eq8} and $q$ vanishes on $\partial\Omega\cup\Sigma$. Finally we obtain
$$\begin{aligned}
\int_\Omega(F_uu'+F_z\cdot\nabla u')\,dx=
&-\int_\Sigma\Big(\frac{\partial u^+}{\partial\nu}
F_z^+\cdot\nu-\frac{\partial u^-}{\partial\nu}F_z^-\cdot\nu\\
&+\frac{\partial u^-}{\partial n}G_u(\nabla q)^-\cdot\nu
-\frac{\partial u^+}{\partial \nu}G_u(\nabla q)^+\cdot\nu\Big)X\nu\,d\HH^1.
\end{aligned}$$
To compute the term $\int_{\Omega}\dive(FX) dx$, we use the divergence theorem thank to the reguality of $\Omega$. Then
$$\int_\Omega\dive(FX)\,dx=\int_{\partial\Omega}FX\nu\,d\HH^{d-1}=0$$
since $X$ is supported in $\Omega$. It follows that
$$\begin{aligned}
\frac{\partial}{\partial\eps}\big|_{\eps=0}
\F(\Sigma_\eps)&=-\lambda\langle H_\Sigma,X\rangle
-\int_\Sigma\Big(\frac{\partial u^+}{\partial\nu}
F_z^+\cdot\nu-\frac{\partial u^-}{\partial\nu}F_z^-\cdot\nu\Big)X\nu\,d\HH^1\\
&-\int_\Sigma\Big(\frac{\partial u^-}{\partial\nu}G_u
(\nabla q)^-\cdot\nu-\frac{\partial u^+}{\partial\nu}G_u(\nabla q)^+\cdot\nu\Big)X\nu\,d\HH^1
\end{aligned}$$
but, by simple computation, we have
$$\begin{aligned}
G_u(\nabla q)^+\cdot\nu&=|\nabla u^+|^{p-2}\nabla q^+\cdot\nu
+(p-2)|\nabla u^+|^{p-4}(\nabla u^+\cdot\nabla q^+)\nabla u^+\cdot\nu\\
&=\Big|\frac{\partial u^+}{\partial\nu}\Big|^{p-2}
\frac{\partial q^+}{\partial\nu}
+(p-2)\Big|\frac{\partial u^+}{\partial\nu}\Big|^{p-4}
\Big(\frac{\partial u^+}{\partial\nu}\frac{\partial q^+}{\partial\nu}\Big)\frac{\partial u^+}{\partial\nu}\\
&=(p-1)\Big|\frac{\partial u^+}{\partial\nu}\Big|^{p-2}\frac{\partial q^+}{\partial\nu}
\end{aligned}$$
and also similarly
$$G_u(\nabla q)^-\cdot\nu=(p-1)\Big|\frac{\partial u^-}{\partial\nu}\Big|^{p-2}\frac{\partial q^-}{\partial\nu}.$$
Combining all the computations together we get
$$\begin{aligned}
\frac{\partial}{\partial\eps}\big|_{\eps=0}&
\F(\Sigma_\eps)=-\lambda\langle H_\Sigma,X\rangle\\
&-\int_\Sigma\Big(\frac{\partial u^+}{\partial\nu}
F_z^+\cdot\nu
-\frac{\partial u^-}{\partial\nu}F_z^-\cdot\nu\Big)X\nu\,d\HH^1\\
&-(p-1)\int_\Sigma\Big(\Big|\frac{\partial u^-}{\partial\nu}\Big|^{p-2}
\frac{\partial u^-}{\partial\nu}\frac{\partial q^-}{\partial\nu}
-\Big|\frac{\partial u^+}{\partial\nu}\Big|^{p-2}
\frac{\partial u^+}{\partial\nu}
\frac{\partial q^+}{\partial\nu}\Big)X\nu\,d\HH^1.
\end{aligned}$$
This equality holds for every vector field $X$, then we derive the following necessary condition of optimality:
$$\begin{aligned}
&-\lambda\langle H_\Sigma,\nu\rangle
-\Big(\frac{\partial u^+}{\partial\nu}
F_z^+\cdot\nu-\frac{\partial u^-}
{\partial\nu}F_z^-\cdot\nu\Big)+\\
&-(p-1)\Big(\Big|\frac{\partial u^-}{\partial\nu}\Big|^{p-2}
\frac{\partial u^-}{\partial\nu}\frac{\partial q^-}{\partial\nu}
-\Big|\frac{\partial u^+}{\partial\nu}\Big|^{p-2}
\frac{\partial u^+}{\partial\nu}\frac{\partial q^+}{\partial\nu}\Big)=0.
\end{aligned}$$
We can rewrite this necessary condition of optimality as
$$\lambda\langle H_\Sigma,\nu\rangle
+\Big(\frac{\partial u}{\partial\nu}F_z\cdot\nu
-(p-1)\Big|\frac{\partial u}{\partial\nu}\Big|^{p-2}
\frac{\partial u}{\partial\nu}
\frac{\partial q}{\partial\nu}\Big)^\pm=0,$$
where the notation $(a)^\pm$ stands for $a^+-a^-$.

We have proved the following result.
\begin{thm}
Let $\Sigma$ be an optimal set in the minimization problem \eqref{eqs1} and $u$ the corresponding solution of the state equation. Assume $d=2$, then $u$ satisfies the necessary condition of optimality:
$$\lambda\langle H_\Sigma,\nu\rangle
+\Big(\frac{\partial u}{\partial\nu}F_z\cdot\nu
-(p-1)\Big|\frac{\partial u}{\partial\nu}\Big|^{p-2}
\frac{\partial u}{\partial\nu}\frac{\partial q}{\partial\nu}
\Big)^\pm=0,$$
where $\nu$ is the unit normal vector of $\Sigma,$ $H_\Sigma$
the generalized mean curvature of $\Sigma$ and $q$
the solution of the adjoint state equation \eqref{eq9}.
\end{thm}

\section{Case of points in $\R^d$, $d>1$}

In the case of points, some extra difficulties
arise because for an equation like \eqref{eq8}
the gradient and the normal are not defined on
points. The equation \eqref{eq8} is not the crucial
point since we are moving an optimal point $x_{0}$ in
the direction of the vector field $X$,
the boundary condition of \eqref{eq8} may be writen
in the general form as $u'(x_{0}) =
\frac{\partial u}{\partial X}(x_{0})$ if $x_0$ is an
optimal point. The main
difficulty is that in the computation we need
an integration by part and since the co-dimention
of the point in $\R^{d} (d>1)$ is greater than $1$,
there is a lack of an integration by part
formulas. The strategy is to study the configurations
which are close to the optimal one and obtain the
optimal configuration as a limit. We consider, for $r$ small and
positive real number, the set $\Sigma_r=\psi(\overline{B_r(x_0)})$
where $B_r(x_0)$ is the ball centered at the point
$x_0$ with radius $r$ and $\psi$ is a smooth
diffeomorphism from $\Omega$ to $\Omega$ such
that $\psi(x_0)=x_0$. The associated state equation is
\begin{equation}\label{eq11}
\begin{cases}
-\dive(|\nabla u|^{p-2}\nabla u)=f&\mbox{in }\Omega\setminus\Sigma_r\\
u=0&\mbox{on }\partial\Omega\cup\Sigma_r.
\end{cases}
\end{equation}
For the functional, we consider
$$\F(\Sigma_r)=\frac{1}{r^{d-1}}
\int_{\Omega\setminus\Sigma_r}F\big(x,u(x),\nabla u(x)\big)\,dx
+\frac{\lambda}{r^{d-1}}\HH^{d-1}(\partial \Sigma_r).$$
The factor $\frac{1}{r^{d-1}}$ is
in order to avoid the functional to degenerate to the
trivial limit functional which vanishes
everywhere. For $r$ small,
$\HH^{d-1}(\partial \Sigma_r) \approx C r^{d-1}$
where $C$ is independent of $r$ and as $r\to 0$ the
solution of the equation \eqref{eq11} converges
strongly in $W^{1, p}_0(\Omega)$ to the solution
of the same equation defined on
$\Omega\setminus\{x_0\}$. If we denote by
$u^{r}$ the solution
of the equation \eqref{eq11} for $r > 0$
 and by $u$ the solution of same equation
for $r = 0$ (solution on
$\Omega \setminus \{x_{0}\}$) we may assume pointwise
convergence of $u^{r}$ and $\nabla u^{r}$
to $u$ and $\nabla u$. Thank to the smoothness of
$F$ we have also the pointwise convergence of
$F(\cdot, u^{r}, \nabla u^{r})$ to
 $F(\cdot, u, \nabla u)$ as $r \to 0$. This gives
heuristically the scaling factor $r^{1-d}$. The same
idea for the case of closed connected sets in the next
section. To simplify the notation in the formulas, we will denote the solution of the equation \eqref{eq11} by $u$ instead of $u^r$. Proceeding as above, that is transforming the domain by $\varphi_\eps$, finding the new state equation and the new functional, and taking the derivative of the functional at $\eps=0$, one gets
$$\frac{\partial}{\partial\eps}\big|_{\eps=0}
\F((\Sigma_r)_\eps)=\frac{1}{r^{d-1}}\int_{\Omega\setminus\Sigma_r}
(F_uu'+F_z\cdot\nabla u'+\dive(FX))\,dx
-\lambda\frac{1}{r^{d-1}}\langle H_{\partial \Sigma_r},X\rangle$$
where $u'$ is the solution of the equation
\begin{equation}\label{eq12}
\begin{cases}
-\dive\big(G_u(\nabla u')\big)=0&\mbox{in }\Omega\setminus\Sigma_r\\
u'=0&\mbox{on }\partial\Omega\\
u'=-\nabla u\cdot X&\mbox{on }\partial \Sigma_r,
\end{cases}
\end{equation}
$G_u(\nabla u')$ is as before and $H_{\partial \Sigma_r}$
is the generalized mean curvature of $\Sigma_r$. Using
the fact that $x_0$ is optimal and $r$ is small enough
(we are in a small neighborhood of the optimal point),
 we obtain $\frac{\partial}{\partial\eps}\big|_{\eps=0}\F((\Sigma_r)_\eps)=o(1)$.

To overcome the problem of $\nabla u'$ as in the
previous case, we introduce the adjoint state equation
\begin{equation}\label{eq13}
\begin{cases}
-\dive\big(G_u(\nabla q)\big)=F_u-\dive(F_z)&\mbox{in }\Omega\setminus\Sigma_r\\
q=0&\mbox{on }\partial\Omega\cup\ \partial \Sigma_r,\\
\end{cases}
\end{equation}
which has to be understood in the distributional sense
$$\int_{\Omega\setminus\Sigma_r}(F_uv-\dive(F_z)v)\,dx
+\int_{\Omega\setminus\Sigma_r}\dive(G_u(\nabla q))v\,dx=0
\qquad\forall v\in C^\infty_c(\Omega\setminus\Sigma_r).$$
In particular
$$\int_{\Omega\setminus\Sigma_r}(F_uu'-\dive(F_z)u')\,dx+
\int_{\Omega\setminus\Sigma_r}\dive(G_u(\nabla q))u'\,dx=0.$$
By integration by parts, the first term of the equation yields
$$\int_{\Omega\setminus\Sigma_r}
(F_uu'-\dive(F_z)u')\,dx=
\int_{\Omega\setminus\Sigma_r}
(F_uu'+F_z\cdot\nabla u')\,dx-
\int_{\partial\Sigma_r}u'F_z\cdot\nu_r\,d\HH^{d-1}$$
where $\nu_r$ is the inward normal of $\Sigma_r$. The
computation is quite similar to the case of
closed connected subset of $\R^2$ and one gets
$$\int_{\Omega\setminus\Sigma_r}\dive(G_u(\nabla q))u'\,dx=
-\int_{\partial\Sigma_r}\Big((p-1)
\Big||\frac{\partial u}{\partial\nu}\Big|^{p-2}
\frac{\partial u}{\partial\nu}\frac{\partial q}{\partial\nu}
\Big)X\nu\,d\HH^{d-1}.$$
Here $\partial \Sigma_{r}$ plays the role of
$\Sigma$ in the two dimensional case. Moreover
all the quantities vanish in the interior side
of $\Sigma_r$ then we are interested only on the
exterior side of $\Sigma_r$. We obtain
$$\int_{\Omega\setminus\Sigma_r}
(F_uu'+F_z\cdot\nabla u')\,dx=-\int_{\partial\Sigma_r}
\Big(\frac{\partial u}{\partial\nu}F_z\cdot\nu
-(p-1)\Big|\frac{\partial u}{\partial\nu}\Big|^{p-1}
\frac{\partial u}{\partial\nu}
\frac{\partial q}{\partial\nu}\Big)X\nu\,d\HH^{d-1}.$$
Using the above calculation, one can
rewrite the derivative of the functional as
$$\begin{aligned}
\frac{\partial}{\partial\eps}\big|_{\eps=0}\F((\Sigma_r)_\eps)=
&-\frac{\lambda}{r^{d-1}}\langle H_{\partial \Sigma_r},X\rangle\\
&-\frac{\lambda}{r^{d-1}}
\int_{\partial\Sigma_r}\Big(\frac{\partial u}
{\partial\nu}F_z\cdot\nu
-(p-1)\Big|\frac{\partial u}{\partial\nu}\Big|^{p-2}
\frac{\partial u}{\partial\nu}
\frac{\partial q}{\partial\nu}\Big)X\nu\,d\HH^{d-1}.
\end{aligned}$$
By a change of variables of type $x=\psi(r,\theta)$, $\theta\in S^{d-1}$ we get
$$\begin{aligned}
-\int_{S^{d-1}}\int_0^r\Big(\frac{\partial u}{\partial\nu}
F_z\cdot\nu-(p-1)\Big|\frac{\partial u}{\partial\nu}\Big|^{p-2}
\frac{\partial u}{\partial\nu}\frac{\partial q}{\partial\nu}\Big)X\nu
J(\theta)\,dr\,d{\theta}
-&\frac{\lambda}{r^{d-1}}\langle H_{\partial\Sigma_r},X\rangle\\
&=o(1).
\end{aligned}$$
In this notation all the terms of the integrand
are evaluated at $\psi(r,\theta)$ and $J(\theta)$
is the Jacobian determinant of the function:
$\theta\mapsto\psi(\theta)$. It remains to study
the limit as $r$ tends to $0$. We do it in the
particular way by letting $\psi(r,\theta)$ go
to $x_0$ in a fixed direction as $r$ goes to $0$. To
express the dependence of the limit on the
direction $\psi(\theta)$, we use the following
notation: $\nu(\psi(r,\theta))\to\nu(\psi(\theta))$
as $r\to0$; the same notation will be also used for
other functions in the integrand. This gives:
\begin{equation}\label{eq14}
\int_{S^{d-1}}\Big(\frac{\partial u}{\partial\nu}
F_z\cdot\nu-(p-1)\Big|\frac{\partial u}{\partial\nu}\Big|^{p-2}
\frac{\partial u}{\partial\nu}
\frac{\partial q}{\partial\nu}\Big)X\nu J\,d\theta=0.
\end{equation}
All the terms in the integrand are evaluated at
$\psi(\theta)$. The quantity
$r^{1-d}\langle H_{\partial \Sigma_r},X\rangle$
goes to zero as $r$ goes to zero. In fact if we
set $\mu_r=r^{1-d}\HH^{d-1}\llcorner \partial \Sigma_r$
this measure belongs to ${\cal M}_{BC}$ and
weakly converges to the Dirac mass $\omega_{d-1}\delta_{x_0}$
concentrated at $x_0$. Since $\mathrm{dim Tan}(\mu_r)=d-1$
for all $r>0$ it follows that $(\mathrm{dim Tan}(\mu_r))
\mu_r=(d-1)\mu_r$ weakly converges to
 $(d-1)\omega_{d-1}\delta_{x_0}$. Therefore
since $\mathrm{dim Tan}(\delta_0)=0<d-1$ we
may apply Theorem \ref{t2} with $g$ the constant
 function $d-1$ to have the weak convergence
of the mean curvature $H_{\mu_{r}}$ to the
mean curvature $H_{\delta_{x_0}}$ which is
identically zero. As a consequence the generalized
 mean curvature $H_{\partial \Sigma_r}$ of $\Sigma_r$ weakly
converges to the generalized mean curvature
$H_{\delta_{x_0}}$ of the point $x_0$. Since
$\partial \Sigma$ is smooth $H_{\partial \Sigma_r}$
coincides with the classical mean curvature but
to avoid confusion with mean curvature of measure
in this paper, we keep the terminology of
general mean curvature. The
equality in \eqref{eq13} holds for every
$X\in C^\infty_c(\Omega)$ and every $\psi$
diffeomorphism. Again it holds true for $X$
constant in the neighborhood of the optimal
point and for all $\psi$ diffeomorphism
satisfying the condition
$$\int_{S^{d-1}}\nu(\psi(\theta))J(\theta)\,d\theta=0.$$
This allows us to write
$$\frac{\partial u}{\partial\nu}F_z\cdot\nu
-(p-1)\Big|\frac{\partial u}{\partial\nu}\Big|^{p-2}
\frac{\partial u}{\partial\nu}\frac{\partial q}{\partial\nu}=const.$$
This expression which is evaluated at $\psi(\theta)$ is constant for all
$\psi$ and $\theta\in S^{d-1}$. This means that it is constant in any direction. Then we have the necessary condition of optimality:
$$\frac{\partial u}{\partial\nu}F_z\cdot\nu
-(p-1)\Big|\frac{\partial u}{\partial\nu}\Big|^{p-2}
\frac{\partial u}{\partial\nu}\frac{\partial q}{\partial\nu}=const.$$
Let us consider a particular case of this problem. We assume $d=2$ and $F=f(x)u$ where $u$ is the solution of the $p$-Laplacian equation. To express the dependence of $u$ on $p$ we denote it by $u_p$ instead of $u$ and the same rule for $q.$ Since $p>2$ we want to study the limit as $p\to2^+$ of the problem. The sequence $u_{p}$ are bounded in $H^1_0(\Omega\setminus\Sigma)$ then up to extracting a subsequence, it converges weakly to some function $u$. It is easy to see that $u$ coincides with the solution of the classical Laplacian that is the solution of equation \eqref{eq1} when $p=2$. From the adjoint state equation \eqref{eq12} we may deduce also that the limit of $q_p$ as $p\to2^+$ coincides with the solution of equation \eqref{eq1} for $p=2$. We may then rewrite the necessary condition of optimality in the form:
$$\Big|\frac{\partial u}{\partial\nu}\Big|=const.$$

The result proved is summarized below.

\begin{thm}
Let $\Sigma$ be an optimal set in the minimization problem \eqref{eqs2}, where $d>1$ and $u$ is the solution of the corresponding state equation. Then $u$ satisfies the necessary condition of optimality:
$$\frac{\partial u}{\partial\nu}F_z\cdot\nu
-(p-1)\Big|\frac{\partial u}{\partial\nu}\Big|^{p-2}
\frac{\partial u}{\partial\nu}
\frac{\partial q}{\partial\nu}= const,$$
where $\nu$ and $q$ are respectively the
limit as $r\to 0$ in a given direction
of the unit normal vector of $\Sigma_r$
and the solution of the adjoint state
equation \eqref{eq13}.
\end{thm}

The case of points in $\R$ is similar
to the case of closed connected subsets in $\R$.

\section{Case of closed connected subsets in $\R^d$ with $d>2$}

Here the strategy is the same. Let $\Sigma$ be the
optimal configuration. We study the configuration
which is close to the optimal one and pass to the
limit. As in the case of points, we consider a
tube $\Sigma_r=\{x\in\R^n\ :\ d(x,\Sigma)\le r\}$. The
associated state equation is
\begin{equation}\label{eq15}
\begin{cases}
-\dive(|\nabla u|^{p-2}\nabla u)=f&\mbox{in }\Omega\setminus\Sigma_r\\
u=0&\mbox{on }\partial\Omega\cup\Sigma_r.\\
\end{cases}
\end{equation}

The procedure is quite similar to the previous
case. The corresponding general functional is
$$\F(\Sigma_r)=\frac{1}{\HH^{d-2}(S^{d-2}_r)}
\int_{\Omega\setminus\Sigma_r}F\big(x,u(x),\nabla u(x)\big)
+\frac{\lambda}{\HH^{d-2}(S^{d-2}_r)}\HH^{d-1}(\partial \Sigma_r),$$
where $S^{d-2}_r$ is a $(d-2)$-dimensional sphere
of radius $r$ and centered on points of
$\Sigma$. From the previous computation,
we deduce the derivative of the functional:
$$\begin{aligned}
\frac{\partial}{\partial\eps}\big|_{\eps=0}&
\F((\Sigma_r)_\eps)=-\frac{\lambda}{\HH^{d-2}(S^{d-2}_r)}
\langle H_{\partial \Sigma_r},X\rangle\\
&-\frac{1}{\HH^{d-2}(S^{d-2}_r)}
\int_{\partial\Sigma_r}
\Big(\frac{\partial u}{\partial\nu}F_z\cdot\nu
-(p-1)\Big|\frac{\partial u}{\partial\nu}\Big|^{p-2}
\frac{\partial u}{\partial\nu}
\frac{\partial q}{\partial\nu}\Big)X\nu\,d\HH^{d-1},
\end{aligned}$$
where $H_{\partial \Sigma_r}$ is the generalized
 mean curvature of $\partial \Sigma_r$. Remark that all
the equations are the same as in the case of
 points in $\R^d$. To pass to the limit,
we argue as in the case of points in $\R^d$. First
we disintegrate the measure $\HH^{d-1}$ and get
$$\begin{aligned}
&\frac{\partial}{\partial\eps}\big|_{\eps=0}
\F((\Sigma_r)_\eps)=-\frac{\lambda}{|\Sigma_r|}
\langle H_{\partial \Sigma_r},X\rangle\\
&-\frac{1}{\HH^{d-2}(S^{d-2}_r)}\int_{\Sigma}\int_{S^{d-2}_r}
\Big(\frac{\partial u}{\partial\nu}F_z\cdot\nu
-(p-1)\Big|\frac{\partial u}{\partial\nu}\Big|^{p-2}
\frac{\partial u}{\partial\nu}
\frac{\partial q}{\partial\nu}\Big)X\nu\,d\HH^{d-2}\,d\HH^1.\\
\end{aligned}$$
The measure $\HH^{d-2}(S^{d-2}_r)^{-1}\HH^{d-2}\llcorner S^{d-2}_r$
converges weakly to $\delta_x$ where $x$ is the center
of the sphere $S^{d-2}_r$ . Due to the hypothesis
made on data of the problem the measure
$$\Big(\frac{\partial u}{\partial\nu}F_z\cdot\nu
-(p-1)\Big|\frac{\partial u}{\partial\nu}\Big|^{p-2}
\frac{\partial u}{\partial\nu}\frac{\partial q}{\partial\nu}\Big)\nu
\HH^{d-2}(S^{d-2}_r)^{-1}
\HH^{d-2}\llcorner S^{d-2}_r$$
weakly converges to the measure
$$\Big(\frac{\partial u}{\partial\nu}F_z\cdot\nu
-(p-1)\Big|\frac{\partial u}{\partial\nu}\Big|^{p-2}
\frac{\partial u}{\partial\nu}
\frac{\partial q}{\partial\nu}\Big)\nu\delta_x$$
for $\HH^1$ a.e. $x\in\Sigma$. The limit here
is computed in a fixed direction
$\theta\in S^{d-2}$. For the curvature part,
notice first that the measure
$\mu_r=\HH^{d-2}(S^{d-2}_r)^{-1}\HH^d\llcorner\partial \Sigma_r\in{\cal M}_{BC}$
weakly converges to the measure
$\HH^1\llcorner\Sigma$ and
$(\mathrm{dim Tan}(\mu_r))\mu_r=d\mu_r$
weakly converges to the measure
 $d \HH^1\llcorner \Sigma$. The
fact that $\mathrm{dim Tan}
\HH^1\llcorner\Sigma=1<d$ allows to apply again
Theorem \ref{t2} to have weak convergence of the
mean curvature of $\mu_r$ to that of
$\HH^1\llcorner\Sigma$ and consequently
the weak convergence of $H_{\partial \Sigma_r}$ to
$H_\Sigma$. Summarizing all computed results we get
$$\begin{aligned}
\frac{\partial}{\partial\eps}\big|_{\eps=0}
\F(\Sigma_\eps)=&-\lambda\langle H_\Sigma,X\rangle\\
&-\int_\Sigma
\Big(\frac{\partial u}{\partial\nu}F_z\cdot\nu
-(p-1)\Big|\frac{\partial u}{\partial\nu}\Big|^{p-2}
\frac{\partial u}{\partial\nu}
\frac{\partial q}{\partial\nu}\Big)X\nu\,d\HH^1=0.
\end{aligned}$$
This relation being true for every vector field $X$,
we get the necessary condition of optimality:
$$\Big(\frac{\partial u}{\partial\nu}F_z\cdot\nu
-(p-1)\Big|\frac{\partial u}{\partial\nu}\Big|^{p-2}
\frac{\partial u}{\partial\nu}
\frac{\partial q}{\partial\nu}\Big)+\lambda
\langle H_\Sigma,\nu\rangle=0.$$

\begin{thm}
Let $\Sigma$ be an optimal set in the minimization
problem \eqref{eqs1}, where $d>2$ and $u$ is the
solution of the corresponding state equation.
Then $u$ satisfies the necessary condition of optimality:
$$\Big(\frac{\partial u}{\partial\nu}F_z\cdot\nu
-(p-1)\Big|\frac{\partial u}{\partial\nu}\Big|^{p-2}
\frac{\partial u}{\partial\nu}\frac{\partial q}{\partial\nu}
\Big)+\lambda \langle H_\Sigma,\nu\rangle=0,$$
where $\nu$ is the unit normal vector of $\Sigma$
in a given direction, $H_\Sigma$ the generalized
 mean curvature of $\Sigma$ and $q$ the limit as
$r\to0$ of the solution of the adjoint state equation.
\end{thm}

Remark that this necessary condition of optimality
depends on the direction $\theta \in S^{d-2}$. Those
directions are contained in the $d-1$ plane which is
 orthogonal to the approximate tangent line to $\Sigma$.

\bigskip
{\bf Acknowledgments.} The work of the second
author has been performed during the Ph. D. studies
 at the Department of Mathematics of University
of Pisa. He also wish to thank the ICTP of Trieste
for a visiting fellowship. The authors wish to thank the anonimous referees for their thorough review and highly appreciate the comments and suggestions, which significantly contributed to improving the quality of the publication.

\bigskip
{\small\noindent
\noindent
Giuseppe Buttazzo:
Dipartimento di Matematica,
Universit\`a di Pisa\\
Largo B. Pontecorvo 5,
56127 Pisa - ITALY\\
{\tt buttazzo@dm.unipi.it}

\bigskip\noindent
Al-hassem Nayam:
International Centre for Theoretical Physics\\
Strada Costiera 11, 34151 Trieste - ITALY\\
{\tt anayam@ictp.it}
}

\end{document}